\documentclass[times]{rncauth}

\usepackage{moreverb}

\usepackage[dvips,colorlinks,bookmarksopen,bookmarksnumbered,citecolor=red,urlcolor=red]{hyperref}
\usepackage{graphicx}
\usepackage{mathrsfs}
\usepackage{amsfonts}
\usepackage{amsmath}
\usepackage{epstopdf}
\newcommand\BibTeX{{\rmfamily B\kern-.05em \textsc{i\kern-.025em b}\kern-.08em
T\kern-.1667em\lower.7ex\hbox{E}\kern-.125emX}}

\def\volumeyear{2023}

 \newtheorem{theorem}{Theorem}

\newtheorem{remark}{Remark}

\begin{document}

\title{Stabilization of  stochastic nonlinear systems
via double-event-triggering
mechanisms and switching controls}

\author{ Xuetao Yang$^{a}$ and Quanxin Zhu$^{*b}$}

\address{$^{a}$ School of Science, Nanjing University of Posts and Telecommunications, Nanjing, 210023, Jiangsu, China \\
         $^{b}$ CHP-LCOCS, School of Mathematics and Statistics,
         Hunan Normal University, Changsha 410081, Hunan, China\\
}


\corraddr{CHP-LCOCS, School of Mathematics and Statistics, Hunan Normal University, Changsha, Hunan 410081, China. \\
E-mail: zqx22@126.com}

\begin{abstract}
In this paper, we concentrate on the exponential stabilization of stochastic nonlinear systems. Different from the single event-triggering mechanism in traditional deterministic/stochastic control systems,
based on two stopping time sequences, we put forward a double-event-triggering
mechanism (DETM) to update control signals and make two different controls  switch in order.
Also,
this novel DETM allows aperiodic time updating
and guarantees a positive lower bound on the inter-event times.
{Together with this DETM, we introduce a switching control law, including a primary control
and a secondary control for a non-switched stochastic system to obtain the exponential stabilization and boundedness results.}
Finally, an illustrative example with simulation figures is given to demonstrate
the obtained results.
\end{abstract}


\maketitle


\vspace{-6pt}

\section{Introduction}
\vspace{-2pt}

In  the passing decades, stability and control problems of stochastic nonlinear systems
 have attracted widespread attention  in
  engineering, neural networks,
physics and so on. {
Nowadays, a large number of control functions are available for stochastic and deterministic systems, like continuous controls, sampled-data controls,
intermittent controls, switching controls and so on (see, e.g.,  \cite{FFMXY20,CH19,M13,GSWW17,GZWPG23,LLLM22,YZ20,LGX19,SYLZ22,LFL07,LL16,LC15}).
In consideration of control cost, sampled-data controls and
intermittent controls are great options. For example,
in \cite{LFL07}, the authors put forward a periodically
intermittent control
$$u(t)= \left\{
\begin{array}{l}
Kx(t),~nT\le t<nT+\tau,\\
0,~~~~~nT+\tau\le t<(n+1)T,
 \end{array}
 \right.
$$
where $K$ is the control gain matrix, $T$ is the control period and $\tau$
is the control width. With the development of intermittent controls (see, e.g., \cite{LL16,LC15}),
periodically intermittent controls are  gradually being replaced by
aperiodic intermittent controls, such as
$$u(t)= \left\{
\begin{array}{l}
Kx(t),~t_i\le t<t_i+\tau,\\
0,~~~~~~~t_i+\tau\le t<t_{i+1},
 \end{array}
 \right.
$$
where $t_{i+1}-t_i\not\equiv T$.
Recently,
Liu et al. in \cite{LYLH21} have applied an important sampling technique called event-triggering mechanism (ETM) to
intermittent controls and proposed
 an aperiodic event-triggered intermittent control as follows:
\begin{eqnarray}
u(t)= \left\{
\begin{array}{l}
Kx(t),~t_i\le t<t_i+\tau_i,\\
0,~~~~~~~t_i+\tau_i\le t<t_{i+1},
 \end{array}
 \right.\label{0}
\end{eqnarray}
where $t_i$ is updated by ETM.
Compared to the traditional time-triggering mechanism, ETM provides aperiodic implementations according to the  triggering conditions. The time will be updated only when triggering conditions are met.
Hence,  ETM can reduce redundant data transmission and efficiently save the network communication resources.
However, in (\ref{0}), $\tau_i$ is not updated by ETM.
Although $\tau_i$ may not be a constant, the choice of $\tau_i$ is difficult under some restrictive
 conditions like (41) and (42) in \cite{LYLH21}, such that in the numerical example, the authors in \cite{LYLH21} took $\tau_i\equiv\tau=1.5$.
Then it is natural to ask whether
  $\tau_i$ can be updated by ETM? In real life,  for example, the light sensing lamps will turn off when they sense strong sunlight and turn on when the sunlight is weak.
Hence, it is more interesting and challenging
to investigate a novel ETM which can generate two time sequences for control updating and switching.

Motivated by the above discussion, we put forward the following double-ETM (DETM):
\begin{eqnarray}
&&s_{i}=\inf\{t> t_{i}\big|~\Psi_1(t,t_i)\ge0\},\nonumber\\
&&t_{i+1}=\inf\{t> s_{i}\big|~\Psi_2(t,s_i)\ge0\},\label{1}
 \end{eqnarray}
where $\Psi_1(t,t_i)$ and $\Psi_2(t,s_i)$ are two  events based on the system state $x(\cdot)$.
DETM (\ref{1}) can determine two execution time sequences $\{t_i\}$, $\{s_i\}$
and they alternate in all stages, i.e., $[0,+\infty)=[t_0,s_0)\cup[s_0,t_1)\cup[t_1,s_1)\ldots$ with $t_0=0$.
As a simple example,  we consider the state function $x=\frac{e\sin3t}{t}+2,~t> 0$ and use the following DETM:
\begin{eqnarray}
&&s_{i}=\inf\{t>t_{i}~\big|~x\le M_2\},\nonumber\\
&&t_{i+1}=\inf\{t>s_{i}~\big|~x\ge M_1\},\label{2}
\end{eqnarray}
where $i=0,1,\ldots$ and $t_0=0$ (see Fig \ref{fig11}).
 \begin{figure}[!h]
\begin{center}
\includegraphics[height=5.7cm]{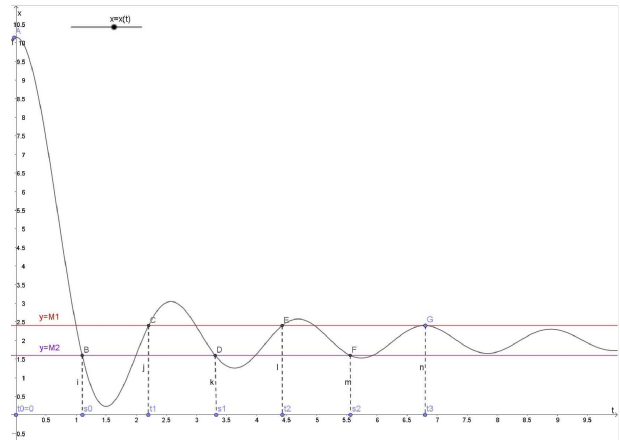}    
\caption{Configuration of DETM.}  
\label{fig11}                                 
\end{center}                                 
\end{figure}

As we know, the classical switching controls  are proposed for switched systems (see,
e.g., \cite{ZF21,NWLXAA19,XPZL20,WNSP17,CZ16,HP21}).
For instance, a switched system $\dot{x}(t)=f_{\sigma(t)}(x_t,u(t))$ is considered in \cite{HP21} and the control $u(t)$ switches according to the
switching signal $\sigma(t)$.
But in this paper,  based on DETM (\ref{1}), we  construct a novel switching control $u(t)$, including a primary control $u_1(t)$ and a secondary control $u_2(t)$
for a non-switched system as follows:
\begin{eqnarray}
u(t)=\left\{
\begin{array}{l}
u_1(t),~~t\in[t_i,s_{i}),\\
u_2(t),~~t\in[s_i,t_{i+1}),
 \end{array}
 \right.  \label{3}
 \end{eqnarray}
where the switching time sequences $\{t_i\}$ and $\{s_i\}$ are both induced by DETM.
Two controls $u_1$ and $u_2$ switch in order according to $\{t_i\}$ and $\{s_i\}$.
Especially, if $u_1\neq0$ and $u_2=0$, then the switching control $u(t)$ degenerates into the intermittent control.
Since most existing results on ETMs  are available for deterministic
control systems (see, e.g.,  \cite{VKPB22,LL22,ZGB22,XTA21,ZLLR22}),
we will face many  challenges and obstacles when applying this switching control and DETM to stochastic systems. This is because
for stochastic  systems, triggering conditions $\Psi_1(t,t_i)$ and $\Psi_2(t,s_i)$ rely on the stochastic process
$x(t)$  such that  $\{t_i\}$ and $\{s_i\}$ are two sequences of stopping times,
 which will cause a lot of difficulties
when considering the stabilization problem and the Zeno phenomenon.
Moreover, although there are a few results on
 stochastic  control systems  (see, e.g., \cite{LD20,ZH21,SNZZL22,QGMY14,L22,WZZ17,Z19}),
 event-triggered controls (ETCs) presented above are all based on a single ETM
and some restrictions are given to avoid the Zeno phenomenon.
Hence, to realize the stabilization of stochastic systems and reduce the restrictive conditions for time regularization at the same time, the design of switching control and DETM requires many stochastic analysis techniques and the
consideration of some new difficulties.}

Consequently, in this paper,
we present a  switching control law based on a novel DETM
 to ensure a desirable stabilization of
stochastic nonlinear systems.  To the best of our knowledge, our paper is the first to put forward this DETM for stochastic control systems and it is novel even in the deterministic control systems.
Compared to the traditional works, the main contributions of the proposed methods are summarized as follows:

{
(1) For a non-switched stochastic system, a novel switching control is established,
which
admits two different controls $u_1$ and $u_2$.  It is more flexible for stabilization of stochastic systems because we can construct $u_1$ and $u_2$,
respectively. Namely, we can choose continuous/discrete controllers  to stabilize
systems or choose two controls with opposite effects, like countermeasure controls.
Especially, if $u_1\neq0$ and $u_2=0$, then this switching control  degenerates into an intermittent control, generalizing those periodic and aperiodic intermittent controls in \cite{LFL07,LL16,LC15,LYLH21}.

(2) Different from the traditional sampled-data controls based on time-triggering mechanisms (see, e.g., \cite{M13,GSWW17,YZ20,LGX19}), in this paper, the switching signal and control updating time are induced by DETM,
which can reduce redundant data transmission and efficiently save the network communication resources.
Also, compared to those event-triggered controls dependent on one   single ETM (see, e.g., \cite{LD20,WZZ17,Z19}), DETM can not only solve the Zeno phenomenon
more effectively, but also
reduce the restrictive conditions for time regularization.
}

The rest of this paper is organized as follows.  In Section 2, we present the notations and
 introduce the model  with an It\^{o} operator.
In Section 3, we propose a detailed design
of DETM
with  execution rules and a novel switching control law based on DETM.
Then, the switching control with different execution  rules is applied to guarantee the  mean square
exponential stabilization and boundedness of stochastic nonlinear systems
in Section 4 and Section 5, respectively.
For illustration of the obtained results,
an example on the pitch motion of a symmetric satellite is provided with simulation figures in Section 6. Finally, the conclusion and future research prospects are
given in the last section.

\vspace{-6pt}

\section{Preliminaries}

\subsection{Notations}
\vspace{-2pt}

In this paper, we use the following notations.
$\mathbf{N}=\{0,1,2,\ldots\}$,
$\mathbf{R}=(-\infty,\infty)$ and $\mathbf{R}^+=[0,\infty)$.
For two positive integers $n$
and $m$, $\mathbf{R}^{n}$ denotes the set of real vectors with the Euclidean norm $|\cdot|$
and $\mathbf{R}^{n\times m}$
denotes the set of real matrices with the induced Euclidean norm $\|\cdot\|$.
Let $\mathbf{C}^i (\mathbf{R}^{n};\mathbf{R}^+)$, $i=0,1,2$ stand for
the space of the continuous functions from $\mathbf{R}^{n}$ to
$\mathbf{R}^+$ with $i$-th order continuous (partial) derivatives.
{ Moreover,
$\alpha\wedge \beta=\min\{\alpha,\beta\}$ and $\alpha\vee \beta=\max\{\alpha,\beta\}$.}

 Let $(\Omega, {\mathcal F}, \{{ \mathcal F}_t\}_{t\geq 0}, P)$ be a complete
probability space with a natural filtration $\{{ \mathcal F}_t\}_{t\geq 0}$
satisfying the usual conditions, i.e.,
$\{{ \mathcal F}_t\}_{t\geq 0}$
 is right continuous and  ${\mathcal F}_0$ contains all $P$-null sets.
$W(t)=(W_1(t),W_2(t),\ldots,W_m(t))^T, ~t\geq0$ is an $m$-dimensional
standard Wiener process
adapted to the filtration $\{{ \mathcal F}_t\}_{t\geq 0}$ on $(\Omega, {\mathcal F}, \{{ \mathcal F}_t\}_{t\geq 0}, P)$.

\subsection{Model description}
\vspace{-0.2cm}
Consider the following  stochastic nonlinear system:
\begin{equation}
\left\{
\begin{array}{l}
dx(t)=f(x(t),u(t))dt+g(x(t))dW(t),~~t>0,\\
x(0)=\phi,
 \end{array}
 \right.  \label{301}
 \end{equation}
where $x\in \mathbf{R}^n$,
$u\in \mathbf{U}\subset \mathbf{R}^l$ is a control input
and $\mathbf{U}$ is a compact set with $0\in \mathbf{U}$.
$f(\cdot,\cdot):\mathbf{R}^n\times\mathbf{R}^l
\rightarrow\mathbf{R}^n$ and
$g(\cdot):\mathbf{R}^n
\rightarrow\mathbf{R}^{n\times m}$ are both continuous functions
satisfying  the local Lipschitz condition
with $f(0,0)=0$ and $g(0)=0$. { Then it follows from \cite{LD20} that system $(\ref{301})$ admits a unique global solution.}

For a function $V\in \mathbf{C}^2(\mathbf{R}^n;\mathbf{R}^+)$,
we present the It\^{o} operator for system (\ref{301}):
\begin{eqnarray*}
\mathcal{L}V(x,u)=\frac{\partial V}{\partial x} f(x,u)+\frac{1}{2}trace(g^T(x)
 \frac{\partial^2 V}{\partial x^2}g(x)),
\end{eqnarray*}
where $\frac{\partial V}{\partial x}=\Big(\frac{\partial V(x)}{\partial x_1},\frac{\partial V(x)}{\partial x_2},\ldots,
\frac{\partial V(x)}{\partial x_n}\Big)$
and $\frac{\partial^2 V}{\partial x^2}=\Big(\frac{\partial^2 V(x)}{\partial x_i\partial x_j}\Big)_{n\times n}$.

\section{Switching control based on DETM}\label{sec-4}

In this section, we introduce a detailed design of DETM with  execution rules
and present a switching control law based on DETM.

\noindent {\bf DETM: the execution rule.} Let $h_1(t),~h_2(t):\mathbf{R}^+\rightarrow \mathbf{R}^+$ be two adjustable functions.
Given a positive constant $\tau$, if the initial value ${|x(0)|^2\le h_2(0)}$,
then we use the following rules:
\begin{eqnarray}
&&t_{i}=\inf\{t>s_{i}~\big|~{|x(t)|^2}\ge h_1(t)\},\nonumber\\
&&s_{i+1}=\inf\{t>t_{i}+\tau~\big|~{|x(t)|^2}\le h_2(t)\},\label{500}
\end{eqnarray}
$i=0,1,\ldots$ with $s_0=0$.


If the initial value ${|x(0)|^2>h_2(0)}$,
then we use the following rules:
\begin{eqnarray}
&&s_{i}=\inf\{t>t_{i}+\tau~\big|~{|x(t)|^2}\le h_2(t)\},\label{5000.1}\\
&&t_{i+1}=\inf\{t>s_{i}~\big|~{|x(t)|^2}\ge h_1(t)\},\label{5000}
\end{eqnarray}
$i=0,1,\ldots$ with $t_0=0$.

Without loss of generality, we consider ${|x(0)|^2>h_2(0)}$
in this paper.

\noindent {\bf Switching control}.
Based on the  time sequences $\{t_i\}_{i\in \mathbf{N}}$ and $\{s_i\}_{i\in \mathbf{N}}$ induced by  execution rules (\ref{5000.1})-(\ref{5000}),
the switching control is designed as:
\begin{eqnarray}
u(t)=\left\{
\begin{array}{l}
u_1(t)=k_1(x(t)),~~t\in[t_i,s_{i}),\\
u_2(t)=k_2(x(t)),~~t\in[s_i,t_{i+1}),
 \end{array}
 \right.  \label{3021}
 \end{eqnarray}
where $k_i:\mathbf{R}^n\rightarrow\mathbf{R}^l$, $i=1,2$
satisfy $k_i(0)=0$ and the local Lipschitz condition.
At each
time instant $t_i$ or $s_i$, the
primary control $u_1(t)$ and the secondary control $u_2(t)$  exchange with each other.
While between consecutive
control updates, the control $u(t)$ stays unchanged.

{\begin{remark}
Similar to  \cite{LD20,FGNZ14,SF16}, we use a time regularization
$\tau$ to avoid the Zeno phenomenon.
But there is no need to add $\tau$ in both
ETM (\ref{5000.1}) and ETM (\ref{5000}). Only one triggering condition needs a positive constant $\tau$.
Namely,
under ETM (\ref{5000.1}),
if execution occurs and the triggering time $t$ is less than $t_i+\tau$, then the primary control $u_1(t)$ will remain working for a fixed period
of time (i.e., $s_i=t_i+\tau$). Otherwise, the triggering time  $s_i=t$.
Hence, the length of each interval $t_{i+1}-t_i=(t_{i+1}-s_i)+(s_i-t_i)$ is clearly
not less than $\tau$, which means that the Zeno phenomenon or even, infinitely fast
execution/sampling cannot occur. Then it is natural to ask which condition to add this  time regularization $\tau$, ETM (\ref{5000.1}) or ETM (\ref{5000})?
A prominent advantage of our DETM is that   the time regularization $\tau$ is added in ETM (\ref{5000.1}) rather than ETM (\ref{5000}), which can avoid the Zeno phenomenon and at the same time,
 reduce the restrictive conditions on  $\tau$ in \cite{LD20,FGNZ14,SF16}.
 The specific reason will be discussed later in Remark \ref{rem4}.
\end{remark}}

\begin{remark}
Compared to the traditional single ETM in \cite{LD20,QGMY14,WZZ17,Z19},
the switching control law (\ref{3021}) based on DETM
admits two different controls $u_1$ and $u_2$ such that it  has more flexibility for control design of stochastic systems. If one of the controls equals 0, then the switching control degenerates into
the aperiodic intermittent control. Also, switching signals $\{t_i\}$
and $\{s_i\}$ in (\ref{3021}) are both induced by ETMs and they are in nature
stochastic (i.e., stopping times), which generate those
intermittent controls in \cite{LFL07,LC15,LYLH21}.
\end{remark}

\begin{figure}[!h]
\begin{center}
\includegraphics[height=3.5cm,width=8cm]{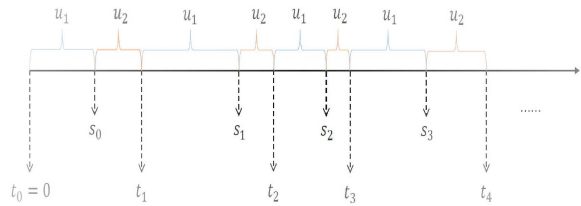}    
\caption{Configuration of switching control based on DETM.}  
\label{fig5.1}                                 
\end{center}                                 
\end{figure}

The configuration
of DETM and switching control
is given by the following procedure (see
 Fig \ref{fig5.1}):

{
Step 1: When $t$ starts from $t_i$, we apply the primary control $u(t)=u_1(t)$;

Step 2: According to  DETM (\ref{5000.1}), we obtain the triggering time
 $t=s_i$ and change the control input, i.e., use the secondary control $u(t)=u_2(t)$ from $t=s_i$ until the next triggering time;

Step 3: Update $t_i$. Based on DETM (\ref{5000}), we obtain $t_{i+1}$
and go into step 1 for $t=t_{i+1}$.
}

%
\section{Exponential stabilization of  stochastic systems }\label{s1}
In this section, we present the  specific expressions for functions $h_1(t)$
 and $h_2(t)$ in the execution rule. Also,  under three key assumptions,
the exponential stabilization is considered for
 system $(\ref{301})$.

In DETM (\ref{5000.1})-(\ref{5000}), take $h_1(t)=M_1e^{-\lambda t}$ and $h_2(t)=M_2e^{-\lambda t}$, i.e.,
\begin{eqnarray}
&&s_{i}=\inf\{t>t_{i}+\tau~\big|~|x(t)|^2\le M_2e^{-\lambda t}\},\label{5000.11}\\
&&t_{i+1}=\inf\{t>s_{i}~\big|~|x(t)|^2\ge M_1e^{-\lambda t}\},\label{5000.12}
\end{eqnarray}
where constants $\lambda>0$, {$0<M_2\le M_1$ and $M_2<|x(0)|^2$.}

\begin{remark}
There are some different constructions of DETM (\ref{5000.11})-(\ref{5000.12}).
For example,
$\lambda$ in ETM (\ref{5000.11}) and ETM (\ref{5000.12}) can be different, i.e.,
\begin{eqnarray*}
&&s_{i}=\inf\{t>t_{i}+\tau~\big|~|x(t)|^2\le M_1e^{-\lambda_2 t}\},\\
&&t_{i+1}=\inf\{t>s_{i}~\big|~|x(t)|^2\ge M_2e^{-\lambda_1 t}\},
\end{eqnarray*}
with $\lambda_1>\lambda_2$.
Also, DETM can be designed as
\begin{eqnarray*}
&&s_{i}=\inf\{t>t_{i}+\tau~\big|~|x(t)|^2\le M_2e^{-\lambda_1 t}\},\\
&&t_{i+1}=\inf\{t>s_{i}~\big|~|x(t)|^2\ge M_1e^{-\lambda_2 t}+\delta(t)\},
\end{eqnarray*}
{ where $\delta(t)$ is an adjustable function satisfying $\lim\limits_{t\rightarrow+\infty}\delta(t)=0$.
Then $\delta(t)$} can be chosen very large in the initial period to reduce  the
working time of the primary control and increase economic efficiency.
\end{remark}

{
 \noindent ${\bf Assumption\  1}$.
There exists a function $V\in \mathbf{C}^2(\mathbf{R}^n;\mathbf{R}^+)$
satisfying
$c_1|x|^2\leq V(x)\leq c_2|x|^2$,
where constants $c_1,c_2>0$.

 \noindent ${\bf Assumption\  2}$.
There exist two
 local Lipschitz mappings $k_1:\mathbf{R}^n\rightarrow \mathbf{U}$
with $k_1(0)=0$
and  $k_2:\mathbf{R}^n\rightarrow \mathbf{U}$
with $k_2(0)=0$.  Moreover,  the following inequalities hold:
\begin{eqnarray}
&&\mathcal{L}V(x,k_1(x))\leq -\mu|x|^2+\bar{M}e^{-\bar{\lambda} t},\label{6000.11}\\
&&\mathcal{L}V(x,k_2(x))\le c_3|x|^2,\label{6000.12}
\end{eqnarray}
where constants  $\mu>0$, $\bar{M}\ge 0$, $\bar{\lambda}<\frac{\mu}{c_2}$ and $c_3\in \mathbf{R}$.

\begin{remark}\label{rem4}
Inequality (\ref{6000.11}) is a classical stable condition while inequality (\ref{6000.12}) is weaker
since $c_3$ can be nonnegative.
Also, it is worth mentioning that, based on  inequality (\ref{6000.11}),  the system will become more stable if time regularization $\tau$ in ETM (\ref{5000.11}) becomes larger. This is because the primary control $k_1$ holds on $[t_i,s_i)$ and $s_i-t_i\ge \tau$.
 Hence, there is no restrictive condition  on the upper bound of $\tau$
 and it is the biggest difference from those traditional ETMs in \cite{LD20,FGNZ14,SF16}.
 If $\tau$ is added in ETM (\ref{5000.12}), then we need similar restrictive conditions
 on $\tau$ in  \cite{LD20,FGNZ14,SF16} because the secondary control $k_2$
  may play a negative role in the stabilization.
 By the way, from the cost of controls,
the less time $k_1$ works, the better.  Then, one can choose a sufficiently small $\tau$  as long as it can avoid  the Zeno phenomenon.
\end{remark}

\begin{remark}
Under the conditions (\ref{6000.11})-(\ref{6000.12}) in Assumption 2,
it can be seen that in DETM (\ref{5000.11})-(\ref{5000.12}), if $M_2$ is chosen smaller, then the primary control $k_1$
works longer and  if $M_1$ is chosen larger, then the secondary control $k_2$
works longer.
\end{remark}

 \noindent ${\bf Assumption\  3}$.
 $\lambda$ in DETM (\ref{5000.11})-(\ref{5000.12}) satisfies
$\lambda<\bar{\lambda}<\frac{\mu}{c_2}$ for $\bar{M}> 0$ and $\lambda<\frac{\mu}{c_2}$ for $\bar{M}= 0$.

\begin{remark}\label{rem2}
Assumption  3 is given to guarantee the existence of
$s_i$ in ETM  (\ref{5000.11}).
Under Assumption  3, if $\bar{M}> 0$, choose
a positive constant $\gamma$ such that $\bar{\lambda}<\gamma<\frac{\mu}{c_2}$. Then for $t_i\le t<s_i$, we have
\begin{eqnarray}
&&\mathbf{E}[e^{\gamma t}V(x(t))-e^{\gamma t_i}V(x(t_i))|{\mathcal F}_{t_i}]\nonumber\\
&&\leq \int_{t_i}^t\mathbf{E}[ \gamma e^{\gamma s}V(x(s))
+e^{\gamma s}\mathcal{L}V(x(s),k_1(x(s)))|{\mathcal F}_{t_i}]ds\nonumber\\
&&\leq c_2 \gamma \int_{t_i}^t e^{\gamma s} \mathbf{E}[|x(s)|^2|{\mathcal F}_{t_i}]ds+ \int_{t_i}^{t}e^{\gamma s}(-\mu\mathbf{E}[|x(s)|^2|{\mathcal F}_{t_i}]+\bar{M}e^{-\bar{\lambda} s})ds\label{41}\\
&&\leq \frac{\bar{M}}{\gamma-\bar{\lambda}}e^{(\gamma-\bar{\lambda})t}.\label{4}
\end{eqnarray}
Taking mathematical expectation on both sides of (\ref{4})
to obtain
\begin{eqnarray*}
e^{\gamma t}\mathbf{E}V(x(t))-e^{\gamma t_i}\mathbf{E}V(x(t_i))
\leq \frac{\bar{M}}{\gamma-\bar{\lambda}}e^{(\gamma-\bar{\lambda})t},
\end{eqnarray*}
which implies
\begin{eqnarray*}
\mathbf{E}|x(t)|^2
&\leq&\frac{c_2 e^{\gamma t_i}}{c_1}\mathbf{E}|x(t_i)|^2 e^{-\gamma t}+
\frac{\bar{M}}{c_1(\gamma-\bar{\lambda})}e^{-\bar{\lambda}t}\\
&\leq&\frac{c_2 e^{\gamma t_i}}{c_1}\mathbf{E}|x(t_i)|^2 e^{-\bar{\lambda} t}+
\frac{\bar{M}}{c_1(\gamma-\bar{\lambda})}e^{-\bar{\lambda}t}\\
&\leq&Le^{-\bar{\lambda} t},
\end{eqnarray*}
 where $L=\max\{\frac{c_2 e^{\gamma t_i}}{c_1}\mathbf{E}|x(t_i)|^2,\frac{\bar{M}}{c_1(\gamma-\bar{\lambda})}\}$.
Then, due to $\lambda<\bar{\lambda}$, one can always find some $t>t_i$ such that
$\mathbf{E}|x(t)|^2\le M_2e^{-\lambda t}$, which shows the existence of $s_i$.
If $\bar{\lambda}\le\lambda$, one may not find $s_0$, let alone $s_i,~i=1,2\ldots$, i.e., the primary
control $k_1$ is applied for $[0,+\infty)$. From the cost of controls,
the less time $k_1$ takes, the better.
By the way, if $\bar{M}=0$, then  $\lambda<\frac{\mu}{c_2}$ and  choose $\lambda<\gamma<\frac{\mu}{c_2}$.
According to the inequality (\ref{41}), we can obtain a similar result.
\end{remark}
}

Now, we are ready to present the main result on stabilization of system
$(\ref{301})$.

\begin{theorem}\label{thm3}
Under { Assumptions 1-3,}
 the solution $x(t)$ to
system $(\ref{301})$ is mean square exponentially stabilized by the switching control (\ref{3021}).
\end{theorem}

\noindent{\bf Proof.}
Based on Remark \ref{rem2}, $s_i$ will definitely be updated after $t_i$.
Hence, there exists some sufficiently large integer $N$ such that $[0,+\infty)=[t_0,s_0)\cup[s_0,t_1)\ldots\cup[s_N,+\infty)$. According to DETM (\ref{5000.11}) and
(\ref{5000.12}), we know that $\mathbf{E}|x(t)|^2\le M_1 e^{-\lambda t}$,
for  $t\in [s_N,+\infty)$. But whether $\mathbf{E}|x(t)|^2\le M_2 e^{-\lambda t}$
for $t\in [s_N,+\infty)$ should be discussed further. We consider the following two cases. \\
Case 1: $\mathbf{E}|x(t)|^2> M_2 e^{-\lambda t},~t\in [s_N,+\infty)$;\\
Case 2: $\mathbf{E}|x(t)|^2\le M_2 e^{-\lambda t},~t\in [s_N,+\infty)$.

For case 1, choose a proper $\gamma$ such that {$\lambda<\gamma<\frac{\mu}{c_2}$.}
Under the Dynkin formula and Assumptions 1-2, we have
\begin{eqnarray}
&&\mathbf{E}e^{\gamma t}V(x(t))-\mathbf{E}e^{\gamma t_0}V(x(t_0))\nonumber\\
&&\leq \mathbf{E}\int_{t_0}^t \gamma e^{\gamma s}V(x(s))
+e^{\gamma s}\mathcal{L}V(x(s),u(s))ds\nonumber\\
&&\leq c_2 \gamma \mathbf{E}\int_{t_0}^t e^{\gamma s} |x(s)|^2ds
+\mathbf{E}\sum_{i=0}^{N} \int_{t_i}^{s_i}e^{\gamma s}\mathcal{L}V(x(s),k_1(x))ds\nonumber\\
&&\quad
+\mathbf{E}\sum_{i=0}^{N-1} \int_{s_i}^{t_{i+1}}e^{\gamma s}\mathcal{L}V(x(s),k_2(x))ds
+\mathbf{E} \int_{s_N}^{t}e^{\gamma s}\mathcal{L}V(x(s),k_2(x))ds\nonumber\\
&&\leq c_2 \gamma \mathbf{E}\int_{t_0}^t e^{\gamma s} |x(s)|^2ds
+\sum_{i=0}^{N-1} \mathbf{E}\int_{s_i}^{t_{i+1}}e^{\gamma s}c_3|x(s)|^2ds\nonumber\\
&&\quad
+\sum_{i=0}^{N}\mathbf{E} \int_{t_i}^{s_i}e^{\gamma s}(-\mu|x(s)|^2+\bar{M}e^{-\bar{\lambda} s})ds
+ \mathbf{E}\int_{s_N}^{t}e^{\gamma s}c_3|x(s)|^2ds.\label{21.2}
\end{eqnarray}

Based on  DETM (\ref{5000.11})-(\ref{5000.12}) and Assumption 3,
take apart the integral interval of $c_2 \gamma \mathbf{E}\int_{t_0}^t e^{\gamma s} |x(s)|^2ds$ and it follows from
(\ref{21.2}) that
\begin{eqnarray}
&&\mathbf{E}e^{\gamma t}V(x(t))-\mathbf{E}e^{\gamma t_0}V(x(t_0))\nonumber\\
&&\leq (c_2 \gamma+c_3) \sum_{i=0}^{N-1} \mathbf{E}\int_{s_i}^{t_{i+1}}e^{\gamma s}|x(s)|^2ds
+(c_2 \gamma-\mu)\sum_{i=0}^{N}\mathbf{E} \int_{t_i}^{s_i}e^{\gamma s}|x(s)|^2ds\nonumber\\
&&\quad+(c_2 \gamma+c_3)\mathbf{E}\int_{s_N}^{t}e^{\gamma s}|x(s)|^2ds
+ \sum_{i=0}^{N} \int_{t_i}^{s_i}e^{\gamma s}\bar{M}e^{-\bar{\lambda} s}ds \label{21.50}
\end{eqnarray}
\begin{eqnarray}
&&\leq (c_2 \gamma+c_3)M_1 \sum_{i=0}^{N-1} \mathbf{E}\int_{s_i}^{t_{i+1}}e^{(\gamma-\lambda)s}ds
+(c_2 \gamma+c_3)M_1\mathbf{E} \int_{s_N}^{t}e^{(\gamma-\lambda)s}ds \nonumber\\
&&\quad+\bar{M} \sum_{i=0}^{N} \mathbf{E}\int_{t_i}^{s_i}e^{(\gamma-\bar{\lambda})s}ds\nonumber\\
&&\leq [(c_2 \gamma+c_3)M_1\vee\bar{M}]\int_{t_0}^{t}e^{(\gamma-\lambda)s}ds\nonumber\\
&&\leq \frac{(c_2 \gamma+c_3)M_1\vee\bar{M}}{\gamma-\lambda}e^{(\gamma-\lambda)t},\label{21.5}
\end{eqnarray}
{which implies
\begin{eqnarray}
\mathbf{E}|x(t)|^2\leq \frac{c_2}{c_1} \mathbf{E}|x(0)|^2 e^{-\gamma t}+\frac{(c_2 \gamma+c_3)M_1\vee\bar{M}}{c_1(\gamma-\lambda)}e^{-\lambda t}\le K e^{-\lambda t},\label{21.6}
\end{eqnarray}
where $K=\max\{\frac{c_2}{c_1} \mathbf{E}|x(0)|^2,\frac{(c_2 \gamma+c_3)M_1\vee\bar{M}}{c_1(\gamma-\lambda)}\}$.}

For case 2, it follows from (\ref{21.50}) that
\begin{eqnarray}
&&\mathbf{E}e^{\gamma t}V(x(t))-\mathbf{E}e^{\gamma t_0}V(x(t_0))\nonumber\\
&&\leq (c_2 \gamma+c_3)M_1 \sum_{i=0}^{N-1} \mathbf{E}\int_{s_i}^{t_{i+1}}e^{(\gamma-\lambda)s}ds
+(c_2 \gamma+c_3)M_2\mathbf{E} \int_{s_N}^{t}e^{(\gamma-\lambda)s}ds \nonumber\\
&&\quad+\bar{M} \sum_{i=0}^{N} \mathbf{E}\int_{t_i}^{s_i}e^{(\gamma-\bar{\lambda})s}ds\nonumber\\
&&\leq [(c_2 \gamma+c_3)M_1\vee\bar{M}]\int_{t_0}^{t}e^{(\gamma-\lambda)s}ds.\nonumber
\end{eqnarray}

Then, it is  the same as (\ref{21.5}) and we can obtain (\ref{21.6}) as well.

This completes the
proof of Theorem \ref{thm3}. $\hfill$ $\Box$

\section{Boundedness of  stochastic systems }
In Section \ref{s1}, we require that $s_{i}-t_i\ge \tau$, i.e., the primary control $u_1$
is active for at least $\tau$-time in the interval $[t_i,t_{i+1})$.
If the minimum inter-event
time $\tau$ is not given in advance, then the following DETM
\begin{eqnarray*}
&&s_{i}=\inf\{t>t_{i}~\big|~|x(t)|^2\le M_2e^{-\lambda t}\},\nonumber\\
&&t_{i+1}=\inf\{t>s_{i}~\big|~|x(t)|^2\ge M_1e^{-\lambda t}\}
\end{eqnarray*}
is hard to avoid the Zeno phenomenon especially for $t\rightarrow+\infty$ even if $M_1>M_2$.

Hence, we put forward an adjustable function $\delta(t)>0$ ($\lim\limits_{t\rightarrow+\infty}\delta(t)\neq 0$) in DETM which can be a  constant function, a  piecewise function
and so on  to avoid the Zeno phenomenon. Namely,
\begin{eqnarray}
&&s_{i}=\inf\{t>t_{i}~\big|~|x(t)|^2\le M_2e^{-\lambda t}\},\nonumber\\
&&t_{i+1}=\inf\{t>s_{i}~\big|~|x(t)|^2\ge M_1e^{-\lambda t}+\delta(t)\},\label{5.1}
\end{eqnarray}
{ where $\lambda>0$, $0<M_2\le M_1+\delta(0)$ and $M_2<|x(0)|^2$.

\begin{theorem}\label{thm3.1}
Under Assumptions 1-3 and DETM (\ref{5.1}),
 the solution $x(t)$ to
system $(\ref{301})$ satisfies: \\
(i) $\mathbf{E}|x(t)|^2\leq L_1$ if $\int_{0}^{t}e^{\lambda (s-t)}\delta(s)ds$
is bounded for $t\ge 0$,
 where $L_1$ is a  positive constant.
Especially, if $\varlimsup\limits_{t\rightarrow+\infty}\delta(t)=\delta$,
 then $\varlimsup\limits_{t\rightarrow+\infty}\mathbf{E}|x(t)|^2\leq(\frac{c_2 }{c_1}+\frac{c_3}{c_1\lambda})\delta$, where $\varlimsup$ denotes the upper limit.\\
(ii) $\mathbf{E}|x(t)|^2
\leq L_2e^{-\lambda t}$ if $\int_{0}^{t}e^{\frac{\mu}{c_2} s}\delta(s)ds$
is bounded for $t\ge 0$,
where $L_2$ is a  positive constant.
\end{theorem}}

\noindent{\bf Proof.}
As the same in Theorem \ref{thm3}, we consider two cases.
Case 1: $\mathbf{E}|x(t)|^2> M_2 e^{-\lambda t},~t\in [s_N,+\infty)$;
Case 2: $\mathbf{E}|x(t)|^2\le M_2 e^{-\lambda t},~t\in [s_N,+\infty)$.

For case 1, let $\lambda<\gamma<\frac{\mu}{c_2}$ and it follows from (\ref{21.50}) that
\begin{eqnarray}
&&\mathbf{E}e^{\gamma t}V(x(t))-\mathbf{E}e^{\gamma t_0}V(x(t_0))\nonumber\\
&&\leq (c_2 \gamma+c_3)M_1 \sum_{i=0}^{N-1} \mathbf{E}\int_{s_i}^{t_{i+1}}e^{(\gamma-\lambda)s}ds
+(c_2 \gamma+c_3)M_1\mathbf{E} \int_{s_N}^{t}e^{(\gamma-\lambda)s}ds \nonumber\\
&&\quad+\bar{M} \sum_{i=0}^{N} \mathbf{E}\int_{t_i}^{s_i}e^{(\gamma-\bar{\lambda})s}ds
+(c_2 \gamma+c_3)\mathbf{E} [\sum_{i=0}^{N-1} \int_{s_i}^{t_{i+1}}
e^{\gamma s}\delta(s)ds
+ \int_{s_N}^{t}e^{\gamma s}\delta(s)ds] \nonumber\\
&&\leq [(c_2 \gamma+c_3)M_1\vee\bar{M}]\int_{t_0}^{t}e^{(\gamma-\lambda)s}ds
+(c_2 \gamma+c_3)\int_{t_0}^{t}e^{\gamma s}\delta(s)ds.\nonumber
\end{eqnarray}

From a similar calculation in the proof of Theorem \ref{thm3},
we have {
\begin{eqnarray*}
\mathbf{E}|x(t)|^2&\leq& Ke^{-\lambda t}
+\frac{c_2 \gamma+c_3}{c_1}e^{-\gamma t}\int_{t_0}^{t}e^{\gamma s}\delta(s)ds\\
&\leq& Ke^{-\lambda t}
+\frac{c_2 \gamma+c_3}{c_1}\int_{t_0}^{t}e^{-\lambda(t-s)}\delta(s)ds,
\end{eqnarray*}
where $K=\max\{\frac{c_2}{c_1} \mathbf{E}|x(0)|^2,\frac{(c_2 \gamma+c_3)M_1\vee\bar{M}}{c_1(\gamma-\lambda)}\}$.}

For case 2, it follows from (\ref{21.50}) again that
\begin{eqnarray}
&&\mathbf{E}e^{\gamma t}V(x(t))-\mathbf{E}e^{\gamma t_0}V(x(t_0))\nonumber\\
&&\leq (c_2 \gamma+c_3)M_1 \sum_{i=0}^{N-1} \mathbf{E}\int_{s_i}^{t_{i+1}}e^{(\gamma-\lambda)s}ds
+(c_2 \gamma+c_3)M_2\mathbf{E} \int_{s_N}^{t}e^{(\gamma-\lambda)s}ds \nonumber\\
&&\quad+\bar{M} \sum_{i=0}^{N} \mathbf{E}\int_{t_i}^{s_i}e^{(\gamma-\bar{\lambda})s}ds
+(c_2 \gamma+c_3)\mathbf{E} \sum_{i=0}^{N-1} \int_{s_i}^{t_{i+1}}
e^{\gamma s}\delta(s)ds \nonumber\\
&&\leq [(c_2 \gamma+c_3)M_1\vee\bar{M}]\int_{t_0}^{t}e^{(\gamma-\lambda)s}ds
+(c_2 \gamma+c_3)\mathbf{E}\int_{t_0}^{t_N}e^{\gamma s}\delta(s)ds.\label{21.08}
\end{eqnarray}

Since the existence of $s_N$, there exists a positive constant $A$ such that $\mathbf{E}\int_{t_0}^{t_N}e^{\gamma s}\delta(s)ds\le A$.
Then, according to (\ref{21.08}) and case 1, we obtain {
\begin{eqnarray*}
\mathbf{E}|x(t)|^2\leq Ke^{-\lambda t} +\frac{A(c_2 \gamma+c_3)}{c_1}e^{-\gamma t}
\leq (K\vee\frac{A(c_2 \gamma+c_3)}{c_1})e^{-\lambda t}.
\end{eqnarray*}
}

This completes the
proof of Theorem \ref{thm3.1}. $\hfill$ $\Box$

{
\begin{remark}
In fact, it is hard to find such a function $\delta(t)$ that satisfies two conditions: $\int_{0}^{t}e^{\frac{\mu}{c_2} s}\delta(s)ds\le$ a positive constant
for all $t\ge 0$
(see  (ii) in Theorem \ref{thm3.1}) and   $\lim\limits_{t\rightarrow+\infty}\delta(t)\neq 0$ (see DETM (\ref{5.1})) at the same time.
Hence, it is difficult to realize the exponential stabilization result of system $(\ref{301})$ . That's why we need
Theorem \ref{thm3} and  DETM  (\ref{5000.11})-(\ref{5000.12}).
\end{remark}}

\section{A numerical example with simulations}
\vspace{-0.2cm}
In this section,
we provide a numerical example
to illustrate  Theorem $\ref{thm3}$.\\
\noindent ${\bf Example}$. Consider the following equation which describes the pitch motion of a symmetric satellite under the influence of aerodynamic torques
and gravity gradient in  \cite{LD20}:
$$\ddot{y}+q_0\dot{y}-l\sin(2y)+q_1(y+q_0\dot{y})\zeta=u,$$
where $y$ is the pitch angle, $u$ is
the control input, { $\zeta$ is Gaussian white noise  satisfying
$\int_{0}^t \zeta(s)ds=W(t),~t\ge 0$,
 $q_1$ is the intensity of the noise, and $q_0, l$ are positive constants.}
Define $x_1=y$, $x_2=\dot{y}$ and  the above system  can be expressed as:
\begin{eqnarray}
 \left(
\begin{array}{l}
dx_1(t)\\
dx_2(t)
 \end{array}
 \right)= \left(
\begin{array}{c}
x_2(t)\\
-q_0x_2(t)+l\sin(2x_1(t))+u(t)
 \end{array}
 \right)dt\nonumber\\
+\left(
\begin{array}{c}
0\\
-q_1(x_1(t)+q_0x_2(t))
 \end{array}
 \right)dW(t).\label{11}
 \end{eqnarray}

For the comparison with \cite{LD20}, we choose the same parameters: $q_0=1, q_1=0.5$ and $l=0.1$.
The switching control is given by {
\begin{eqnarray*}
u(t)=\left\{
\begin{array}{l}
k_1(x)=-x_1(t),~t\in[t_i,s_{i}),\\
k_2(x)=0,~~~~~~~~~~t\in[s_i,t_{i+1}).
 \end{array}
 \right.
 \end{eqnarray*}}

 The switching times $\{t_i\}$ and $\{s_i\}$ based on DETM will be defined later.
Choose the Lyapunov function $V(x)=\frac{47}{36}x_1^2+\frac{10}{9}x_1x_2+
x_2^2-\frac{1}{5}\sin^2(x_1)$ for the system. From a simple calculation, we obtain $c_1=0.47,~c_2=1.73$ in Assumption 1. Also,
\begin{eqnarray*}
{\mathcal{L} V(x,k_1(x))}
=-\frac{31}{36}x_1^2+\frac{1}{9}x_1\sin(2x_1)-\frac{23}{36}x_2^2
\leq-\frac{23}{36}(x_1^2+x_2^2)\leq -\frac{1}{3}|x|^2,
\end{eqnarray*}
and
\begin{eqnarray*}
{\mathcal{L} V(x,k_2(x))}
=\frac{1}{4}x_1^2+\frac{1}{9}x_1\sin(2x_1)-\frac{23}{36}x_2^2+2x_1x_2
\leq\frac{53}{36}x_1^2+\frac{13}{36}x_2^2\leq 2|x|^2.
\end{eqnarray*}

Then, it is easy to see that
$\mu=\frac{1}{3}$, $\bar{M}=0$ and $c_3=2$ { in Assumption 2.}

Let
the initial value $x(0)=[-2,3]^T$.
In DETM (\ref{5000.11})-(\ref{5000.12}),
choose $M_1=0.8,M_2=0.5,\lambda=0.1<\frac{\mu}{c_2},~\tau=0.5$.
Then, { Assumption 3 holds.}
By Theorem \ref{thm3},
we
conclude that system $(\ref{11})$ is exponentially stable in mean square.

\begin{figure}[!h]
\begin{center}
\includegraphics[height=4.5cm,width=7cm]{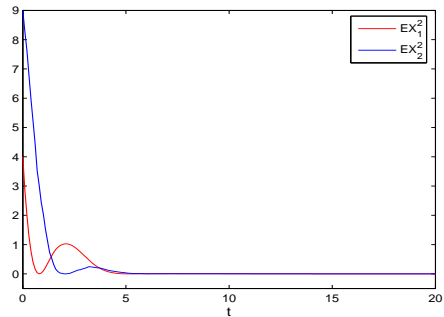}    
\caption{Second moment of the solution to Eq.(\ref{11}) under DETM.}  
\label{fig1}                                 
\end{center}                                 
\end{figure}

\begin{figure}[!h]
\begin{center}
\includegraphics[height=4.5cm,width=7cm]{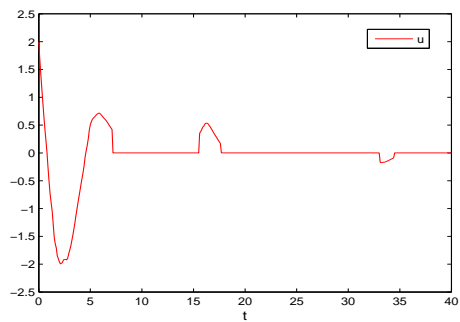}    
\caption{Response of switching control $u(t)$.}  
\label{fig2}                                 
\end{center}                                 
\end{figure}

\begin{figure}[!h]
\begin{center}
\includegraphics[height=4.5cm,width=7cm]{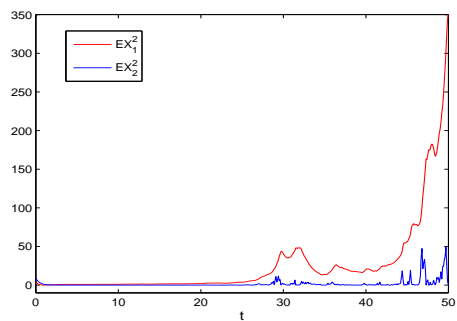}
\caption{Second moment of the solution to Eq.(\ref{11}) without control.}  
\label{fig3}                                 
\end{center}                                 
\end{figure}

\begin{figure}[!h]
\begin{center}
\includegraphics[height=4.5cm,width=7cm]{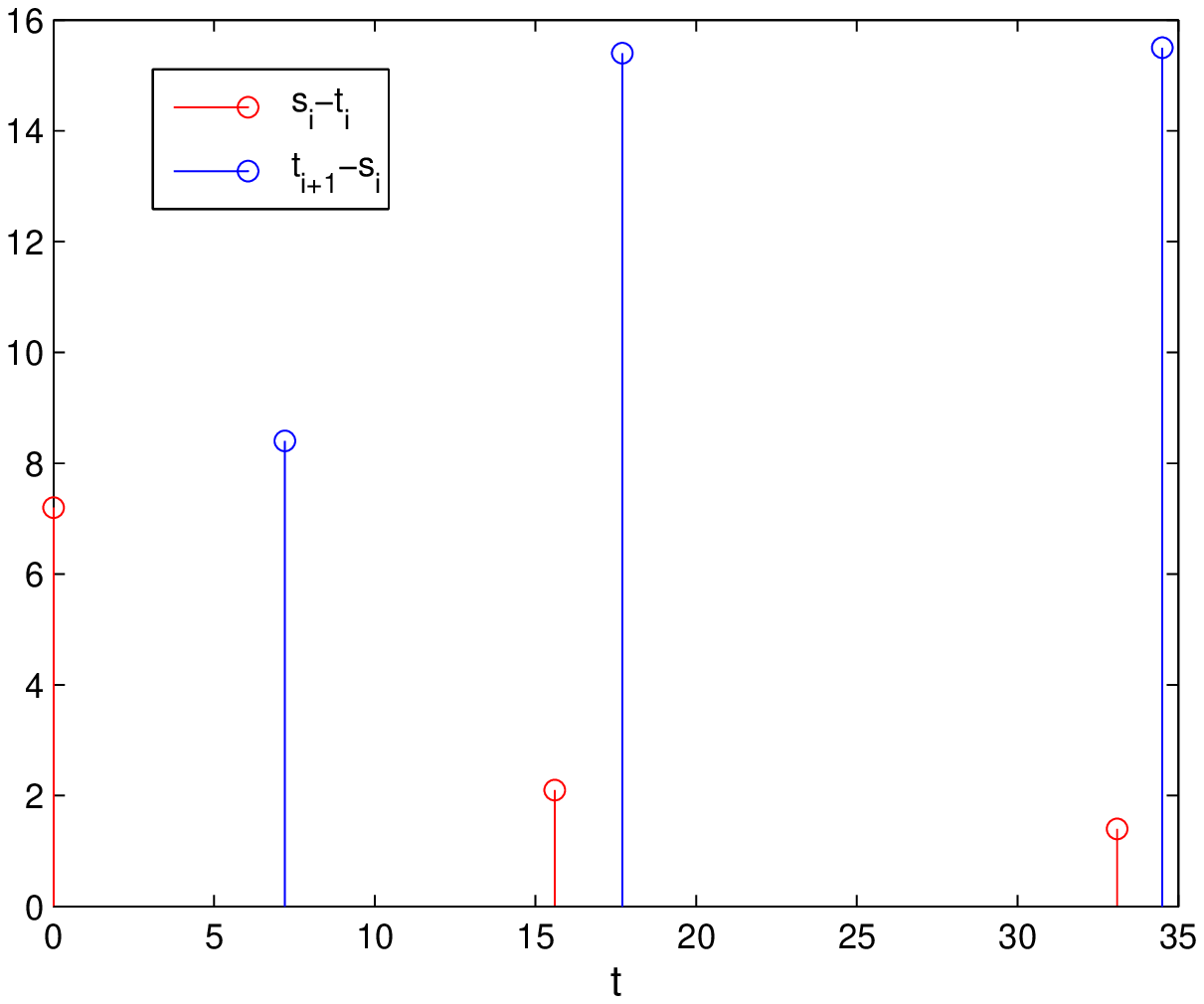}
\caption{ Inter-event times $s_{i}-t_i$ and $t_{i+1}-s_i$.}  
\label{fig4}                                 
\end{center}                                 
\end{figure}

The simulations are given as follows.
Using Matlab with step=0.1 and $[0,T]=[0,50]$, 100
sample path trajectories are simulated.
We draw the following
four figures. Fig \ref{fig1} and Fig \ref{fig3} show the time
evolutions of  mean
square of the solutions to system $(\ref{11})$  and the corresponding
uncontrolled system, respectively. It is obvious
that system $(\ref{11})$  is mean square exponentially stable and the corresponding
uncontrolled system is unstable.
 Fig \ref{fig2}  draws
the behaviour of the switching control $u(t)$.
{ It can be seen that the switching signals $\{t_i\}$ and $\{s_i\}$ are stochastic
and both induced by ETMs. Hence, this switching control $u(t)$ generalizes those periodic and aperiodic intermittent controls in \cite{LFL07,LL16,LC15,LYLH21}.}
Fig \ref{fig4} plots the inter-event
times $s_{i}-t_i$ and $t_{i+1}-s_i$, { which shows both triggering times and duration  of the primary control are very small. Compared to ETCs given in the example in \cite{LD20}, whatever the static ETC or dynamic ETC,}
the number of triggering times is much more than that induced by DETM
(see Fig. (c) and Table I in \cite{LD20}). Hence, the switching control dependent on DETM
can reduce the
communication times greatly.

\section{Conclusions and future research prospects}

To obtain the stabilization and boundedness of   stochastic nonlinear systems, we propose a switching control law based on a novel DETM. The advantage of
 DETM is that it can generate two  stopping time series
to realize the control updating and switching.
Also, in DETM, the minimum inter-event time
is guaranteed to be bounded by a positive constant to avoid the Zeno phenomenon.
The switching control law based on DETM
admits two different controls and extends the traditional
 aperiodic intermittent control.
However, in this paper, the secondary control has little effect on stabilization
because we choose a continuous  primary control.
In the future work, we will consider the primary control to be zero-order-hold.
Then, the effectiveness of the  secondary control will be highlighted.

\bibliographystyle{plain}

\begin{thebibliography}{22}

\bibitem{FFMXY20}
C. Fei, W. Fei, X. Mao, D. Xia and L. Yan.
 Stabilization of highly nonlinear hybrid systems by feedback control based on discrete-time state observations.
\textit{ IEEE Trans. Automat.
Control}, 2020, 65: 2899-2912.



\bibitem{CH19}
A. Cetinkaya and T. Hayakawa.
 A sampled-data approach to pyragas-type delayed feedback stabilization of periodic orbits.
\textit{ IEEE Trans. Automat.
Control}, 2019, 64: 3748-3755.

\bibitem{M13}
X. Mao.
 Stabilization of continuous-time hybrid stochastic differential
equations by discrete-time feedback control.
\textit{ Automatica}, 2013, 49: 3677-3681.

\bibitem{GSWW17}
Y.  Gao, X.  Sun, C. Wen and W. Wang.
 Estimation of sampling
period for stochastic nonlinear sampled-data systems with emulated controllers.
\textit{ IEEE Trans. Automat.
Control}, 2017,
62: 4713-4718.

\bibitem{GZWPG23}
X.  Guo, D.  Zhang, J.  Wang, J.  Park and L. Guo.
 Observer-based event-triggered composite anti-disturbance control for multi-agent systems under multiple disturbances and stochastic FDIAs.
\textit{ IEEE Trans. Autom. Sci. Eng.},  2023, 20:
528-540.

\bibitem{LLLM22}
X. Li, W. Liu, Q. Luo and X. Mao.
Stabilisation in distribution of hybrid stochastic differential equations by feedback control based on discrete-time state observations.
\textit{ Automatica}, 2022,  140,  110210.

\bibitem{YZ20} X. Yang and Q. Zhu.
 Stabilization of stochastic retarded systems based on sampled-data feedback control.
\textit{ IEEE Trans. Syst., Man, Cybern., Syst.}, 2021, 51: 5895-5904.



\bibitem{LGX19}
S. Li, J. Guo and Z. Xiang.
 Global stabilization of a class of switched
nonlinear systems under sampled-data control.
\textit{ IEEE Trans. Syst.,
Man, Cybern., Syst.}, 2019,  49: 1912-1919.



\bibitem{SYLZ22} J.  Sun, J. Yang, S.  Li and Z.  Zeng.
Predictor-based periodic event-triggered control for dual-rate networked control systems with disturbances.
\textit{  IEEE Trans.  Cybern.},  2022, 52:
8179-8190.





\bibitem{LFL07}
C. Li, G. Feng and X. Liao.
Stabilization of nonlinear systems via periodically intermittent control.
\textit{  IEEE Trans. Circuits Syst. II, Exp. Briefs}, 2007,
54: 1019-1023.



{
\bibitem{LL16}  Y. Li and C. Li.
Complete synchronization of delayed chaotic
neural networks by intermittent control with two switches in a
control period. \textit{ Neurocomputing}, 2016, 173: 1341-1347.



\bibitem{LC15} X.  Liu and T.  Chen.
Synchronization of complex networks
via aperiodically intermittent pinning control. \textit{ IEEE Trans. Automat.
Control}, 2015, 60: 3316-3321.
}

 \bibitem{LYLH21}
B. Liu, M. Yang, T. Liu and D. Hill.
 Stabilization to exponential input-to-state stability
via aperiodic intermittent control.
\textit{  IEEE Trans. Automat. Control}, 2021,
66: 2913-2919.







\bibitem{ZF21}
Y. Zhao and  J. Fu.
 $H_\infty$ composite anti-bump switching
control for switched systems. \textit{  IEEE Trans. Syst.,
Man, Cybern., Syst.}, Early Access, 2021,
1-9.

\bibitem{NWLXAA19} B. Niu, D. Wang, H. Li, X. Xie, N.  Alotaibi and F.  Alsaadi.
A novel neural-network-based adaptive control scheme for outputconstrained
stochastic switched nonlinear systems.
\textit{  IEEE Trans. Syst.,
Man, Cybern., Syst.}, 2019, 49: 418-432.

\bibitem{XPZL20}
X.  Xiao, J.  Park, L. Zhou and G.  Lu.
 Event-triggered control of discrete-time switched linear systems with
network transmission delays.
\textit{  Automatica}, 2020,
111, 108585.


\bibitem{WNSP17}
 J. Wen, S. Nguang, P. Shi and L. Peng.
Finite-time stabilization of Markovian jump delay systems-a switching control approach.
\textit{ Int. J.
Robust Nonlinear Control}, 2017, 27: 298-318.


\bibitem{CZ16}
Y. Chen and W. Zheng.
Stability analysis and control for switched stochastic delayed systems.
\textit{ Int. J.
Robust Nonlinear Control}, 2016, 26: 303-328.

\bibitem{HP21}
I. Haidar and   P. Pepe.
 Lyapunov-Krasovskii charaterization of the input-to-state stability for switching retarded systems.
\textit{  SIAM J. Control Optim.}, 2021, 59: 2997-3016.


\bibitem{VKPB22}
C. Viel, M. Kieffer, H. Piet-Lahanier and S. Bertrand.
 Distributed event-triggered formation control for multi-agent systems in presence of packet losses.
\textit{  Automatica}, 2022, 141,  110215.







\bibitem{LL22}
X. Li and P. Li.
Input-to-state stability of nonlinear systems: event-triggered impulsive control.
\textit{  IEEE Trans. Automat. Control}, 2022, 67: 1460-1465.

\bibitem{ZGB22}
K. Zhang, B. Gharesifard and E. Braverman.
 Event-triggered control for nonlinear time-delay systems.
\textit{  IEEE Trans. Automat. Control}, 2022, 67: 1031-1037.


\bibitem{XTA21}
X. Xu, A. Tahir and B. A\c{c}{\i}kme\c{s}e.
Periodic event-triggered control for incrementally quadratic nonlinear systems.
\textit{ Int. J.
Robust Nonlinear Control}, 2021, 31: 5261-5280.



\bibitem{ZLLR22}
H. Zhao, W. Li, Z. Li and Y. Ren.
Dynamic event-triggered control for Markovian jump linear systems with time-varying delay.
 \textit{ Int. J.
Robust Nonlinear Control}, 2022, 32: 2359-2379.


\bibitem{LD20}
S. Luo and F. Deng.
 On event-triggered control of nonlinear stochastic systems.
\textit{  IEEE Trans. Automat. Control}, 2020, 65: 369-375.


\bibitem{ZH21}
Q. Zhu and T. Huang.
 $H_\infty$ control of stochastic networked control systems with time-varying delays: The event-triggered sampling case.
 \textit{ Int. J.
Robust Nonlinear Control}, 2021, 31: 9767-9781.



\bibitem{SNZZL22}
W. Su, B. Niu, Y. Zhang, J. Zhang and Y. Li.
Event-triggered adaptive decentralized control of stochastic nonlinear systems with strong interconnections and time-varying full-state constraints.
 \textit{ Int. J.
Robust Nonlinear Control}, 2022, 32: 5200-5225.

\bibitem{QGMY14} D. Quevedo, V. Gupta, W. Ma and S. Yuksel.
 Stochastic stability of
event-triggered anytime control.
\textit{  IEEE Trans. Automat. Control}, 2014, 59: 3373-3379.

{ \bibitem{L22}
S. Luo.
Stability and $\mathcal{L}_2$-gain analysis of linear periodic
event-triggered systems with large delays.
\textit{ IEEE Trans. Circuits Syst. II, Exp. Briefs}, Early Access, 2022, 1-5.}

\bibitem{WZZ17} Y. Wang, W. Zheng and H. Zhang.
 Dynamic event-based control of
nonlinear stochastic systems.
\textit{  IEEE Trans. Automat. Control}, 2017, 62: 6544-6551.



\bibitem{Z19}
Q. Zhu.
Stabilization of stochastic nonlinear delay systems with exogenous disturbances and the event-triggered feedback control.
\textit{  IEEE Trans. Automat. Control}, 2019,
 64: 3764-3771.






\bibitem{FGNZ14}
F. Forni, S. Galeani, D. Ne\v{s}i\'{c} and L. Zaccarian.
 Event-triggered transmission
for linear control over communication channels.
\textit{  Automatica}, 2014,
50: 490-498.


\bibitem{SF16}
 A. Selivanov and E. Fridman.
 Event-triggered $H^\infty$ control: A switching
approach.
\textit{  IEEE Trans. Automat. Control}, 2016,
61: 3221-3226.







\end{thebibliography}

\end{document}